\documentclass[11pt]{article}
\usepackage[left=1in,right=1in,top=1in,bottom=1in]{geometry}

\usepackage{amsmath,amssymb}
\usepackage{theorem}
\usepackage{color}
\usepackage{mathrsfs}

\theorembodyfont{\itshape}
\theoremstyle{plain}
\newtheorem{theorem}{Theorem}[section]

\newtheorem{lemma}[theorem]{Lemma}
\newtheorem{fact}[theorem]{Fact}
\newtheorem{corollary}[theorem]{Corollary}
\newtheorem{proposition}[theorem]{Proposition}
\theorembodyfont{\rmfamily}
\newtheorem{example}[theorem]{Example}
\newtheorem{remark}[theorem]{Remark}
\newtheorem{question}[theorem]{Question}
\newtheorem{problem}[theorem]{Problem}
\newtheorem{definition}[theorem]{Definition}

\newtheorem{claim}{Claim}

\newcommand{\id}{\mbox{{\rm id}}}

\def\Z{\mathbb{Z}}
\def\Q{\mathbb{Q}}
\def\R{\mathbb{R}}

\begin{document}

\title{Topological groups all continuous automorphisms of which are open}
\author{Vitalij A.~Chatyrko and Dmitri B.~Shakhmatov}

\date{}

\maketitle

\begin{center}
{Dedicated to Professor Michael Megrelishvili on the occasion of his 60th anniversary}
\end{center}

\smallskip

\begin{abstract}
A topological space is reversible if  each continuous bijection of it onto itself is open. We introduce an analogue of this notion in the category of topological groups:
A topological group $G$ is {\em $g$-reversible\/} if every continuous automorphism of $G$ (=continuous isomorphism of $G$ onto itself) is open.  The class of $g$-reversible groups contains Polish groups, locally compact $\sigma$-compact groups, minimal groups, abelian groups with the Bohr topology, and reversible topological groups. We prove that subgroups of 
$\R^n$ 
are $g$-reversible, for every positive integer $n$. An example of a compact (so reversible) metric abelian group having a countable dense non-$g$-reversible subgroup is given. 
We also highlight the differences between reversible spaces and $g$-reversible topological groups. Many open problems are scattered throughout the paper.
\end{abstract} 


\smallskip
{\it 2000 AMS (MOS) Subj. Class.:} Primary: 22A05; Secondary: 20K30, 22B05, 22D05, 46A30, 54A10, 54D45, 54H11
\medskip

\bigskip
{\em All topological groups considered in this paper are supposed to be Hausdorff.\/}

The symbol $\mathbb{R}$ denotes the additive group of real numbers in its usual topology. By $\mathbb{Q}$ and $\mathbb{Z}$ we denote subgroups of $\R$ consisting of rational numbers and integer numbers, respectively (in their subspace topology).
The quotient group $\mathbb{T}=\mathbb{R}/\mathbb{Z}$ is called the {\em torus group\/}.

We refer the reader to \cite{AT}, \cite{HM}, \cite{M} or \cite{D}
for necessary information on topological groups, and to \cite{E} for undefined topological notions.

\section{Elementary facts}

This preliminary section contains notations and elementary facts which shall be needed in Section~\ref{sec:2}. The proofs are included only for the convenience of the readers.

\begin{definition}
\label{def:two:topologies}
Let $X$ be a set and $f:X\to X$ be a bijection. For every topology $\tau$ on $X$ we consider two topologies
\begin{equation}
\label{eq:1}
\tau_f^\rightarrow=\{f(U): U\in\tau\}
\
\text{ and }
\
\tau_f^\leftarrow=\{f^{-1}(U): U\in\tau\}
\end{equation}
on $X$.
\end{definition}
\begin{remark}
\label{rem:1.2}
{\em In the notation of Definition \ref{def:two:topologies},
both maps 
$f:(X,\tau_f^\leftarrow)\to (X,\tau)$ and 
$f:(X,\tau)\to (X, \tau_f^\rightarrow)$ are homeomorphisms.}
Indeed, 
since $f$ is a bijection, from the second equation in \eqref{eq:1} we get
\begin{equation}
f(\tau_f^\leftarrow)=\{f(V):V\in \tau_f^\leftarrow\}
=
\{f(f^{-1}(U)):U\in\tau\}
=
\{U:U\in\tau\}=\tau.
\end{equation}
Furthermore,
\begin{equation}
f(\tau)=\{f(U):U\in\tau\}=\tau_f^\rightarrow
\end{equation}
by the first equation in \eqref{eq:1}.
\end{remark}

\begin{fact}
\label{fact:1.3}
For a bijection $f:X\to X$ of a set $X$ and a topology $\tau$ on $X$, the following conditions are equivalent:
\begin{itemize}
\item[(i)] the map $f:(X,\tau)\to (X,\tau)$ is continuous;
\item[(ii)] $\tau_f^\leftarrow\subseteq \tau$;
\item[(iii)] $\tau\subseteq \tau_f^\rightarrow$.
\end{itemize}
\end{fact}
Proof.
It follows from the second equation in \eqref{eq:1} that (i) and (ii) are equivalent.

(ii)$\to$(iii)
Since $\tau_f^\leftarrow\subseteq \tau$, we have 
$\tau=f(\tau_f^\leftarrow)\subseteq f(\tau)=\tau_f^\rightarrow$ by \eqref{eq:1}.

(iii)$\to$(ii)
Since 
$\tau\subseteq \tau_f^\rightarrow$ and $f$ is a bijection, we have
$$
\tau_f^\leftarrow=\{f^{-1}(U):U\in\tau\}=f^{-1}(\tau)\subseteq f^{-1}(\tau_f^\rightarrow)=
\{f^{-1}(f(U)):U\in\tau\}=\{U: U\in\tau\}=\tau
$$
by \eqref{eq:1}.
$\Box$

\begin{fact}
\label{fact:1.4}
For a bijection $f:X\to X$ of a set $X$ and a topology $\tau$ on $X$, the following conditions are equivalent:
\begin{itemize}
\item[(i)] the map $f:(X,\tau)\to (X,\tau)$ is open;
\item[(ii)] $\tau_f^\rightarrow\subseteq \tau$;
\item[(iii)] $\tau\subseteq \tau_f^\leftarrow$.
\end{itemize}
\end{fact}
Proof. Item (i) is equivalent to the continuity of the inverse map $g=f^{-1}$ of $f$. Applying Fact  \ref{fact:1.3} to the bijection $g$,
we conclude that 
\begin{equation}
\label{eq:4}
(i)\leftrightarrow \tau_g^\leftarrow\subseteq \tau
\leftrightarrow \tau\subseteq \tau_g^\rightarrow.
\end{equation}
It remains to note that
$$
\tau_g^\leftarrow=\{g^{-1}(U): U\in\tau\}=\{f(U):U\in\tau\}=\tau_f^\rightarrow
$$ 
and
$$
\tau_g^\rightarrow=\{g(U): U\in\tau\}=\{f^{-1}(U):U\in\tau\}=\tau_f^\leftarrow.
$$
Combining this with \eqref{eq:4} and Definition \ref{def:two:topologies}, we obtain the conclusion of our claim.
$\Box$

\section{Introduction} 
\label{sec:2}

In \cite{RW} Rajagopalan and  Wilansky gave
the following definition:
\begin{definition}
\label{def:reversible}
A topological space $X$ is called {\em reversible\/} if
each continuous bijection of $X$ onto itself is open (or equivalently, a homeomorphism). 
\end{definition}

Note that the openness of a bijection $f$ of a topological space $X$ on itself is equivalent to continuity of its inverse map $f^{-1}$. 
This observation is probably at the origin of the term ``reversible'',
as the continuity of every bijection $f$ of a reversible space automatically implies the continuity of its inverse $f^{-1}$; that is, the continuity of a bijection becomes ``reversible''.    
   
Rajagopalan and  Wilansky proved the following fact in \cite{RW}:
\begin{fact}
\label{characterization:of:reversibility}
For a topological space $(X,\tau)$, the following conditions are equivalent:
\begin{itemize}
\item[(i)] $(X,\tau)$ is not reversible;
\item[(ii)] there exists a topology $\tau_w$ on $X$ such that  
$\tau_w \subsetneq \tau$ and the space $(X, \tau_w)$ is homeomorphic to  $(X, \tau)$;
\item[(iii)] there exists a topology $\tau_s$ on $X$ such that  
$\tau \subsetneq \tau_s$ and the space $(X, \tau_s)$ is homeomorphic to  $(X, \tau)$.
\end{itemize}
\end{fact}

\medskip
Proof.
By Definition \ref{def:reversible},
item (i) is equivalent to the existence of a continuous bijection 
$f:(X,\tau)\to (X,\tau)$ such that $f$ is not open. 
By Fact \ref{fact:1.3}, the continuity of $f$ is equivalent to each of the two inclusions
$\tau_f^\leftarrow\subseteq \tau$ and $\tau\subseteq \tau_f^\rightarrow$.
By Fact \ref{fact:1.4}, the non-openness of $f$ is equivalent 
to each of the two non-inclusions
$\tau_f^\rightarrow\not\subseteq \tau$ and $\tau\not\subseteq \tau_f^\leftarrow$.
This means that item (i) is equivalent to each of the two formulae
 $\tau_f^\leftarrow\subsetneq \tau$ and $\tau\subsetneq \tau_f^\rightarrow$. By Remark \ref{rem:1.2},
both of the spaces $(X,\tau_f^\leftarrow)$ and $(X,\tau_f^\rightarrow)$ are homeomorphic to $(X,\tau)$.
So one can take $\tau_f^\leftarrow$ as $\tau_w$ and $\tau_f^\rightarrow$ as $\tau_s$.
$\Box$

\medskip
  
Let us recall some known notions from the theory of (topological) groups.  
  
\begin{definition}
\begin{itemize}
\item[(i)] 
Given two groups $G_1$ and $G_2$, a map $h: G_1\to G_2$ is said to be an {\em 
isomorphism\/} provided that $h$ is a bijection which is also a group homomorphism. 
\item[(ii)]
An isomorphism of a group $G$ onto itself is called an {\em automorphism\/} of $G$.

\item[(iii)] A map $h: (G_1,\tau_1) \to (G_2,\tau_2)$ between topological groups $(G_1,\tau_1)$ and $(G_2,\tau_2)$
is called a 
{\em topological isomorphism\/} if $h$ is both an 
isomorphism and a homeomorphism. 
\item[(iv)]
Topological groups $(G_1, \tau_1)$ and $(G_2, \tau_2)$ are said to be {\em topologically isomorphic\/}
provided that there exists a topological isomorphism $h: (G_1,\tau_1) \to (G_2,\tau_2)$ 
between them.
\end{itemize}
\end{definition}
\vskip 0.3 cm

The following notion is a natural analog of the notion of reversibility 
for topological groups:

\begin{definition}
\label{def:g-reversible}
We say that a topological group $G$ is {\em $g$-reversible\/} (an abbreviation for {\em group reversible\/}) if every continuous automorphism of $G$ is open (and thus, a topological isomorphism of $G$ onto itself).
\end{definition}

As seen from this definition, $g$-reversibility of a topological group can be viewed as some sort of an ``open mapping property''.

\begin{remark}
\label{inverse:map:remark}
Clearly, a topological group $G$ is $g$-reversible if and only if every continuous automorphism $f$ of $G$ has continuous inverse $f^{-1}$. 
\end{remark}

\begin{remark}
Let $G$ be a topological group. As usual, $\mathrm{Aut}(G)$ denotes the group of all automorphisms of $G$. Denote by 
$\mathrm{Aut}_c(G)$ the monoid of all continuous automorphisms
of $G$ and by $\mathrm{Aut}_t(G)$ the subgroup of $\mathrm{Aut}(G)$ consisting of all topological isomorphisms of $G$ onto itself. Clearly,
$\mathrm{Aut}_t(G)\subseteq \mathrm{Aut}_c(G)\subseteq \mathrm{Aut}(G)$. Note that the equality $\mathrm{Aut}_t(G)= \mathrm{Aut}_c(G)$ holds precisely
when $G$ is $g$-reversible. It follows that $G$ is $g$-reversible 
if and only if the monoid $\mathrm{Aut}_c(G)$ is a group.
\end{remark}

Due to an obvious similarity between Definitions \ref{def:reversible} and \ref{def:g-reversible}, it comes as no surprise that 
the following analog of 
Fact
\ref{characterization:of:reversibility} holds in the category of topological groups.

\begin{proposition}
\label{two:topologies}
For a topological group $(G,\tau)$, the following conditions are equivalent:
\begin{itemize}
\item[(i)] $(G,\tau)$ is not $g$-reversible;
\item[(ii)] there exists a group topology $\tau_w$ on $G$ such that  
$\tau_w \subsetneq \tau$ and the topological group $(G, \tau_w)$ is topologically isomorphic to  $(G, \tau)$;
\item[(iii)] there exists a topology $\tau_s$ on $G$ such that  
$\tau \subsetneq \tau_s$ and the topological group $(G, \tau_s)$ is topologically isomorphic to  $(G, \tau)$.
\end{itemize}
\end{proposition}
Proof.
By Definition \ref{def:g-reversible},
item (i) is equivalent to the existence of a continuous automorphism
$f:(G,\tau)\to (G,\tau)$ such that $f$ is not open. 
By Fact \ref{fact:1.3}, the continuity of $f$ is equivalent to each of the two inclusions
$\tau_f^\leftarrow\subseteq \tau$ and $\tau\subseteq \tau_f^\rightarrow$.
By Fact \ref{fact:1.4}, the non-openness of $f$ is equivalent 
to each of the two non-inclusions
$\tau_f^\rightarrow\not\subseteq \tau$ and $\tau\not\subseteq \tau_f^\leftarrow$.
This means that item (i) is equivalent to each of the two formulae
 $\tau_f^\leftarrow\subsetneq \tau$ and $\tau\subsetneq \tau_f^\rightarrow$. Since $f$ is an automorphism of $G$, so is its inverse $f^{-1}$. Therefore, both $\tau_f^\leftarrow$ and $\tau_f^\rightarrow$ are (Hausdorff) group topologies on $G$ and 
 both maps 
$f:(G,\tau_f^\leftarrow)\to (G,\tau)$ and 
$f:(G,\tau)\to (G, \tau_f^\rightarrow)$ are topological isomorphisms.
So one can take $\tau_f^\leftarrow$ as $\tau_w$ and $\tau_f^\rightarrow$ as $\tau_s$.
$\Box$

\medskip
Since each automorphism of a group $G$ is a also
a bijection of $G$ onto itself, from Definitions \ref{def:reversible}
and \ref{def:g-reversible}, one obtains the following 

\begin{proposition}\label{reversible_g_reversible}
Each reversible topological group is $g$-reversible.
\end{proposition}

Let us give a simple example highlighting the difference 
between the notions of reversibility and $g$-reversibility, as well as showing that 
the converse of Proposition~\ref{reversible_g_reversible} is not valid.

\begin{example}\label{the group of rational numbers} 
Let $\mathbb{Q}$ be the group of rational numbers in its 
usual topology inherited from the real line $\R$.
It is known (see \cite[Example 3]{RW}) that {\em the topological space $\mathbb Q$ is not reversible\/}. 

Let us show that {\em the topological group
$\mathbb Q$ is $g$-reversible\/}.  In fact, consider any continuous automorphism $f: \mathbb Q \to \mathbb Q$. Let $r = f(1)$. 
Since $f$ is an automorphism, $r=f(1)\not =f(0)=0$.
It is easy to see that $f(q) = r \cdot q$ for each $q \in \mathbb Q$. 
Note that the mapping $h(q) : \mathbb Q \to \mathbb Q$ defined by $h(q) = r^{-1} \cdot q$ for each $q \in \mathbb Q$, is continuous and coincides with the inverse of $f$. 
By 
Remark~\ref{inverse:map:remark},
the topological group $\mathbb Q$ is $g$-reversible. 
\end{example}

Next, we recall basic facts about reversible spaces.

\begin{fact}
\label{reversible:fact}
\begin{itemize}
\item[(i)]
discrete spaces are reversible;
\item[(ii)]
compact Hausdorff spaces are reversible;
\item[(iii)] locally Euclidean spaces are reversible;
\item[(iv)] neither the space of rational numbers $\Q$ nor the space of irrational numbers $\mathbb{P}$ (considered with the subspace topology inherited from the reals $\R$ with its usual topology) are reversible;
\item[(v)] the topological union $N_{\aleph_0} = D_{\aleph_0} \oplus cD_{\aleph_0}$ of the countably infinite discrete space $D_{\aleph_0}$ 
with its one point compactification $cD_{\aleph_0}$ is not reversible;
\item[(vi)] for connected reversible spaces $X$ and $Y$, the topological union $X \oplus Y$  is always reversible;
\item[(vii)] the topological product $X \times Y$ of reversible spaces $X$ and $Y$ need not be reversible, even when both $X$ and $Y$ are connected;
\item[(viii)] if a space $X$ is not reversible, then
the topological union $X \oplus Y$ and the topological product $X \times Y$ are not reversible, for any space $Y$.
\end{itemize}
\end{fact}

Item  (vii) of this fact is shown in \cite{ChK}, while the rest is taken from \cite{RW}.

\begin{remark}
\label{rem:2.10}
It follows from Proposition \ref{reversible_g_reversible} and items (i)--(iii) of Fact \ref{reversible:fact} that
the following topological groups are $g$-reversible:
\begin{itemize}
\item
discrete groups, 
\item
compact groups,
\item
the additive group $\R^n$ of an $n$-dimensional vector space
taken with its Euclidean topology (for a positive integer $n$),
\item
topological products $D \times \mathbb R^n$ of a discrete group $D$ and $\R^n$, for some positive integer $n$.
\end{itemize}
\end{remark}

\bigskip

In this article 
we thoroughly investigate the new notion of $g$-reversibility 
and compare it with the old notion of reversibility.
In particular, we present numerous examples of topological groups which are not reversible as topological spaces but $g$-reversible as topological groups.

An overview of the paper follows.
In Section~\ref{sec:3} we show that many classical topological groups are $g$-reversible. The list includes Polish groups, separable metric groups with the automatic continuity property, minimal groups and abelian groups with their Bohr topology.
Furthermore,
$\sigma$-compact locally compact groups are $g$-reversible as well (Theorem~\ref{omega_narrow_theorem}); in particular, locally compact connected groups are $g$-reversible (Corollary~\ref{connected:locally:compact}). Nevertheless, (either locally connected or zero-dimensional) locally compact metric groups need not be $g$-reversible (Corollary~\ref{corollary:4.4}).
In Section~\ref{sec:5} we show that the Euclidean groups $\R^n$ behave much better in this respect; indeed, {\em every\/} subgroup of $\R^n$ is 
$g$-reversible (Theorem~\ref{theorem_all_subgroups_R_n}).

In Section~\ref{sec:6}, we give an example
of a countable dense subgroup of the Hilbert space which is not $g$-reversible (Example ~\ref{separable_metrizable_countable_non-g-reversible}) and we show that a countable dense subgroup of a compact abelian group need not be $g$-reversible either 
(Theorem~\ref{precompact:countable:non-sm}).
In Section~\ref{sec:7} we study topological groups every subgroup of which is $g$-reversible in the subgroup topology.
In Section~\ref{sec7:b} we deduce from a theorem of Megrelishvili that {\em every\/} topological group is a group retract of some $g$-reversible group. From another theorem of Megrelishvili, we deduce that a closed subgroup of a locally compact metric $g$-reversible group need not be $g$-reversible (Example~\ref{loc:comp:ex:of:closed:subgroup}).

In Section~\ref{sec:8} we study a question of preservation of $g$-reversibility under taking Tychonoff products. 
In Section~\ref{sec:9} we give four examples contrasting reversibility with $g$-reversibility in particular topological groups. 

We formulate many open problems which are scattered throughout the whole paper. The last Section~\ref{sec:10} has some
remarks and additional open problems.
In particular, in Remark~\ref{Q:remarks}
we notice that the space $\mathbb Q$ of rational numbers admits two group structures such that one topological group is $g$-reversible and the other one is not.

\section{Many classical topological groups are $g$-reversible}
\label{sec:3}

In this section we shall demonstrate that many classical topological groups are $g$-reversible.

A topological space is {\em \v{C}ech-complete\/} if (and only if) it is a $G_\delta$ subset of its Stone-\v{C}ech compactification.

A topological group $G$ is called {\em $\omega$-narrow\/} 
provided that for every open neighbourhood $U$ of the identity 
of $G$ one can find an at most countable 
subset $S$ of $G$ such that $G=SU$ \cite{Guran}.
A topological group is $\omega$-narrow if and only if
it is topologically 
isomorphic to a subgroup of a suitable product of separable metric groups \cite{Guran}.

\begin{theorem}
\label{omega-narrow:Cech-complete}
Every  $\omega$-narrow \v{C}ech-complete
group is $g$-reversible.
\end{theorem}
\medskip
Proof. 
This  follows from the fact that every continuous surjective homomorphism between $\omega$-narrow \v{C}ech-complete
groups is open \cite[Corollary 4.3.33]{AT}.
$\Box$

\medskip
Let us recall that {\it a Polish space\/} is a separable space the topology of which can be generated by a complete metric.
In particular, Polish spaces are separable metrizable spaces, so have a countable base. 

Since separable metric groups are $\omega$-narrow and Polish spaces are \v{C}ech-complete~\cite[Theorem 4.3.26]{E}, from Theorem \ref{omega-narrow:Cech-complete} we get the following corollary.

\begin{corollary}\label{Polish_groups} Every Polish 
group is $g$-reversible.  
\end{corollary}

Recall that a topological group $G$ has {\it the automatic continuity property\/} if every group homomorphism from $G$ to  a separable 
metric 
group is continuous  
(cf. \cite{Ros}).                                                                                                                     
In  \cite{Ros, S} one finds  numerous examples of  separable
metric
groups with the automatic continuity property.

The following proposition easily follows 
from  
                                                                                                                     Remark \ref{inverse:map:remark}.
\begin{proposition}
A
separable 
metric
group with the automatic continuity property is $g$-reversible. 
\end{proposition}

\medskip
Recall that a topological group $G$ is {\em minimal\/} if every continuous isomorphism from $G$ onto a (Hausdorff) group $H$ is open. 
Combining this with Definition~\ref{def:g-reversible}, we obtain the following

\begin{theorem}
\label{minimal:theorem}
Every minimal 
group is $g$-reversible. 
\end{theorem}

\begin{remark}
It follows from Definition~\ref{def:g-reversible} that $g$-reversibility of a topological group $G$ can be viewed as a weaker version of minimality obtained by restricting the group $H$ in the definition of a minimal group to be the group $G$ itself. Since only the continuous isomorphisms 
of $G$ onto {\em itself\/} are required to be open in Definition~\ref{def:g-reversible},  $g$-minimal groups might perhaps be called ``{\em self\/}-minimal'', with Theorem~\ref{minimal:theorem} stating that ``minimal $\to$ self-minimal''.
\end{remark}

Minimal groups are extensively studied in topological group theory; see the survey \cite{DM} and references therein.

Let us recall that a Hausdorff abelian group $G$ is called {\it precompact\/} if $G$ is a subgroup of a Hausdorff compact group, or equivalently, for each open neighbourhood $U$ of zero of $G$ there exists a finite subset $F$ of $G$ such that $U+F = G$. 

Every abelian group $G$ has the strongest precompact group topology on $G$ called its {\em Bohr topology\/}. 
This topology is the initial topology for the family of all homomorphisms from $G$ to the torus group $\mathbb{T}=\mathbb{R}/\mathbb{Z}$;
such homomorphisms of $G$ are called characters of $G$. Therefore, the Bohr topology of $G$ is the coarsest topology on $G$ making all characters of $G$ continuous (cf. \cite[Theorem 1]{M} or \cite{D}).

\begin{theorem}\label{Bohr_topology}
Every 
abelian group with the Bohr topology is  $g$-reversible.
\end{theorem}
Proof.
Indeed, let $G$ be an abelian group and $\tau$ be the Bohr topology of $G$. Assume that $(G, \tau)$ is not $g$-reversible.
By Proposition \ref{two:topologies}, 
there exists a group topology $\tau_s$ on the group $G$ such that $\tau \subsetneq \tau_s$ and the topological group $(G, \tau_s)$ is topologically isomorphic to $(G, \tau)$. 
Since $(G,\tau)$ is precompact and precompactness is preserved by topological isomorphisms, $(G, \tau_s)$ is precompact.
Now the inclusion $\tau \subsetneq \tau_s$ contradicts 
the fact that $\tau$ is the strongest precompact group topology on $G$. 
$\Box$
                                                                                                                       
                                                                                                                                                                                                                                                                                                                                                                       \section{$g$-reversibility in locally compact groups}
                                                                                                                                                                                                                                                                                                                                                                        \label{sec:4}
                                                                                                                                                                                                                                                                                                                                                                       
\begin{theorem}\label{omega_narrow_theorem}
Every  $\sigma$-compact locally compact group is $g$-reversible.
\end{theorem}
Proof. The statement directly follows from Definition 
 \ref{def:g-reversible} and an old result of Pontryagin 
 saying that every continuous homomorphism  of a $\sigma$-compact topological group  onto a locally compact topological group is open (cf. \cite[Theorem 3.1.27]{AT}). 
$\Box$

\medskip

Local compactness in Theorem~\ref{omega_narrow_theorem} is essential and cannot be omitted,
even in the presence of metrizability;
this will be shown in 
Theorem \ref{precompact:countable:non-sm}.

Since connected locally compact groups are $\sigma$-compact
\cite[Corollary 3.1.5]{AT},
 we get the following

\begin{corollary}
\label{connected:locally:compact}
A connected locally compact group is $g$-reversible.
\end{corollary}

In Question V from \cite[Section 9]{RW}, Rajagopalan and Wilansky
ask whether every connected locally compact group must be reversible.
In view of Proposition \ref{reversible_g_reversible}, Corollary \ref{connected:locally:compact}
provides a positive answer to a weaker ($g$-reversible) version of this question.

                                                                                                                                                                                                                                                                                                                                                                        Connectedness in Corollary \ref{connected:locally:compact} cannot be replaced with local connectedness; see Corollary \ref{corollary:4.4}(ii) below.
                                                                                                                                                                                                                                                                                                                                                                        
                                                                                                                                                                                                                                                                                                                                                                        As was mentioned in Remark in \cite{Ra}, first examples of continuous automorphisms $f$ of some locally compact (abelian) groups which are not open, have been given by Robertson in his Ph.D. thesis \cite{Rob}. Another example, based on an idea of Hofmann, was given in \cite{Ra} and \cite[Example 7]{RW}. We include here a slightly more general construction.

\begin{example}
\label{example:RW}
{\em Define $Z_{+}=\{z\in \Z:z>0\}$ and $Z_{-}=\Z\setminus Z_{+}$. Let $K$ be an infinite compact group. Consider $D=K^{Z_{-}}$ with the discrete topology and $C=K^{Z_{+}}$ with the Tychonoff product topology.
Then the product $G=D\times C$ is not $g$-reversible.\/} 
Indeed, it
is easy to check that the right shift map
$f:G\to G$ defined by $f((x_n))=(x_{n-1})$ for $(x_n)\in G$,
is a continuous automorphism of $G$.
The converse $g=f^{-1}$ of $f$ is the left shift map defined by
$g((x_n))=(x_{n+1})$ for $(x_n)\in G$. 
Note that $L=\prod_{n\le 0} \{e\}\times K\times \prod_{n\ge 2} \{e\}$ is an infinite compact subset of $G$, while its image 
$g(L)= \prod_{n< 0} \{e\}\times K\times \prod_{n\ge 1} \{e\}$
under $g$
is an infinite discrete subset of $G$, so it cannot be compact.
This shows that $g$ is discontinuous.
By Remark \ref{inverse:map:remark}, $G$ is not reversible.
\end{example}

\begin{corollary}
\label{g-nonreversible_product}
\begin{itemize}
\item[(i)]
The product $D\times C$ of a discrete abelian group $D$ of order 2 and cardinality
continuum and a zero-dimensional compact metric abelian group $C$ of order 2 need not be $g$-reversible. 
\item[(ii)] The product $D\times C$ of a discrete abelian group $D$ of cardinality
continuum and a connected, locally connected compact metric abelian group $C$ need not be $g$-reversible. 
\end{itemize} 
\end{corollary}

Proof. To prove item (i), we apply Example \ref{example:RW}
to $K=\{0,1\}^\Z$, and to prove item (ii), we apply the same example to $K=\mathbb{T}^\Z$, where $\mathbb{T}=\R/\Z$ is the torus group.
$\Box$

\begin{corollary}
\label{corollary:4.4}
\begin{itemize}
\item[(i)]
A zero-dimensional locally compact metric abelian group (of order $2$) need not be $g$-reversible.
\item[(ii)]
A locally connected locally compact metric abelian group need not be $g$-reversible.
\end{itemize}
\end{corollary}

This corollary shows that $\sigma$-compactness in Theorem~\ref{omega_narrow_theorem} 
and connectedness in Corollary \ref{connected:locally:compact}
are essential conditions and cannot be omitted, even in the presence of metrizability.

\begin{problem}
Describe $g$-reversible locally compact (abelian) groups.
\end{problem}

\section{$g$-reversibility in subgroups of Euclidean groups $\R^n$} 
\label{sec:5}

We have seen in the last section that locally compact (abelian) groups need not be $g$-reversible. The goal of this section is to show that the Euclidean groups behave much better in this respect; namely, we prove in Theorem~\ref{theorem_all_subgroups_R_n}  that {\em every\/} subgroup of the Euclidean group is 
$g$-reversible.

\begin{fact}
{\rm (\cite[Proposition 7.5(ii)]{HM})}
\label{fact:5.1}
Let $E_1$, $E_2$ be topological vector spaces over $\R$. 
Then  every continuous group homomorphism
$h: (E_1,+)\to (E_2,+)$ between their additive groups
is a linear map.
\end{fact}

\begin{fact}
\label{fact:dikranjan}
{\rm (\cite[Theorem 4.37]{D})}
Let $n$ be a non-negative integer.
Every
closed subgroup of  $\mathbb R^n$  is topologically isomorphic to a product  $\Z^l \times \R^m$ for some non-negative integers $l$ and $m$ such that $l+m \leq n$.
\end{fact}

The following fact is a classical result from linear algebra.

\begin{fact}
\label{fact:linear:surjections}
Let $n$ be a non-negative integer.
Every surjective linear map $f:\R^n\to\R^n$ has an inverse 
$f^{-1}$ which is also a linear map.
\end{fact}

\begin{corollary}
\label{monomorphism}
Let $n$ be a non-negative integer.
Every surjective homomorphism
$\varphi: \Z^n\to \Z^n$
is a monomorphism (and thus, an automorphism of $\Z^n)$.
\end{corollary}

Proof. 
Note that there exists a linear map $f:\R^n\to \R^n$ such that 
$f\restriction_{\Z^n}=\varphi$. Since the image 
$f(\R^n)$ of $\R^n$ under $f$ contains $f(\Z^n)=\varphi(\Z^n)=\Z^n$ by our assumption, 
the dimension 
of the image $f(\R^n)$ is equal to $n$.
Since $f(\R^n)$ is a linear subspace of $\R^n$ having the same 
dimension as $\R^n$, it follows that $f(\R^n)=\R^n$;
that is, $f$ is surjective.
Since $f:\R^n\to\R^n$ is a surjective linear map, it
must have trivial kernel by Fact \ref{fact:linear:surjections}.
We conclude that $f$ is a monomorphism, and so is 
its restriction $f\restriction_{\Z^n}=\varphi$ to $\Z^n$.
$\Box$

\begin{lemma}
\label{extention:lemma}
Let $H=C\times D$ be a product of a connected group $C$ and 
a discrete group $D$. Suppose that $h:H\to H$ is a continuous homomorphism such that $h(H)$ is dense in $H$.
Then there exists a surjective homomorphism $\varphi: D\to D$ such that 
$p\circ h=\varphi\circ p$, where $p: H=C\times D\to D$ is the projection on the second coordinate. 
\end{lemma}

Proof. Let $d\in D$ be arbitrary. Since both $h$ and $p$ are continuous maps, so is their composition $p\circ h: H\to D$.
Since $C$ is  connected by our assumption, so is $C\times\{d\}$.
Therefore, its image $p\circ h(C\times\{d\})$ under the continuous map $p\circ h$ is a connected subset of $D$.
Since $D$ is discrete, $p\circ h(C\times\{d\})$ must be a singleton. Therefore, there exists a unique element $\varphi(d)\in D$ such that 
\begin{equation}
\label{eq:5:k}
p\circ h(C\times\{d\})=\{\varphi(d)\};
\text{ that is, }
h(C\times\{d\})\subseteq C\times \{\varphi(d)\}.
\end{equation}
This defines a map $\varphi: D\to D$.

Let $(c,d)\in H$. Then $c\in C$ and $d\in D$, so applying \eqref{eq:5:k}, we get
$p\circ h(c,d)=\varphi(d)=\varphi(p(c,d))=\varphi\circ p(c,d)$.
This proves that $p\circ h=\varphi\circ p$.

To show that $\varphi$ is surjective, fix an arbitrary $y\in D$.
Since $D$ is discrete, $C\times\{y\}$ is a non-empty open subset of $H$. Since $h(H)$ is dense in $H$, we can find $g\in H$ such that
$h(g)\in C\times\{y\}$. There exist $c\in C$ and $d\in D$ such that
$g=(c,d)$. Then $g\in C\times \{d\}$ and 
$h(g)\in h(C\times\{d\})\subseteq C\times\{\varphi(d)\}$ by 
\eqref{eq:5:k}.
This means that $\varphi(d)=y$.

It remains only to check that $\varphi$ is a homomorphism.
Let $d_1,d_2\in D$ be arbitrary. Since both $h$ and $p$ are homomorphisms, from \eqref{eq:5:k} we obtain that
$$
\varphi(d_1d_2^{-1})=p\circ h(e, d_1d_2^{-1})=p\circ h(e,d_1) \cdot (p\circ h(e,d_2))^{-1}=\varphi(d_1)\cdot \varphi(d_2)^{-1}.
$$
This proves that $\varphi$ is a homomorphism.
$\Box$

\medskip
The following theorem is a main result in this section.

\begin{theorem}\label{theorem_all_subgroups_R_n} 
Let $n$ be a non-negative integer.
Then
every subgroup of $\mathbb R^n$ 
is $g$-reversible.
\end{theorem}
Proof. 
Let $G$ be a subgroup of $\mathbb R^n$. 
The closure $H$ of $G$ in $\R^n$ is Raikov complete, being a closed subgroup of the Raikov complete group $\R^n$.
Therefore, $H$ coincides with the Raikov completion of $G$.

Let $f : G \to G$ be a continuous automorphism of $G$.
Since every continuous homomorphism admits a unique extension to a continuous homomorphism between Raikov completions of its domain and image, respectively, the exists a continuous homomorphism $h: H\to H$ extending $f$. 
Note that
$G=f(G)=h(G)\subseteq h(H)\subseteq H$, and since 
$G$ is dense in $H$, we conclude that $h(H)$ is dense in $H$.

Since $H$ is a closed subgroup of $\R^n$, 
we can use 
Fact \ref{fact:dikranjan}
to fix 
non-negative integers $l$ and $m$ such that $H$ is topologically isomorphic to the product $C\times D$, where $C=\R^m$ and $D=\Z^l$. Without loss of generality, we shall assume that $H=C\times D$. Since $C$ is connected and $D$ is discrete, all the assumptions of Lemma \ref{extention:lemma} are satisfied.

Applying Lemma \ref{extention:lemma}, we can fix a homomorphism  $\varphi$ as in the conclusion of this lemma.
Since $D=\Z^l$ and $\varphi$ is surjective, $\varphi$ is a monomorphism by Corollary \ref{monomorphism}.

Define 
\begin{equation}
\label{eq:6:h}
C_0=C\times\{0\}
\text{ and }
G_0=G\cap C_0.
\end{equation}
Since $\varphi(0)=0$, we have 
\begin{equation}
\label{eq:7d}
h(C_0)=h(C\times\{0\})\subseteq C\times \{\varphi(0)\}
=C\times\{0\}=C_0.
\end{equation}
\begin{claim}
\label{claim:1}
$f_0=f\restriction_{G_0}:G_0\to G_0$ is a surjection of $G_0$ onto itself.
\end{claim}

Proof.
Since $h$ extends $f$,
from \eqref{eq:6:h} and \eqref{eq:7d} we get
$f(G_0)=h(G_0)\subseteq h(C_0)\subseteq C_0$
and 
$f(G_0)\subseteq f(G)=G$,
which implies
$f(G_0)\subseteq G\cap C_0=G_0$.
Therefore,
$f_0=f\restriction_{G_0}:G_0\to G_0$ is a map of  $G_0$ into itself.

Let us verify that $f_0: G_0\to G_0$ is surjective.
Indeed, let $g_0\in G_0$ be chosen arbitrarily.
Since $f$ is an automorphism of $G$ and $g_0\in G_0\subseteq G$, there exists $g_1\in G$ such that $f(g_1)=g_0$.

Let $p$ be the projection defined in Lemma \ref{extention:lemma}.
Since $g_0\in G_0\subseteq C_0=C\times\{0\}$ by \eqref{eq:6:h},
the property of $\varphi$ implies that
$$
\varphi(p(g_1))=\varphi\circ p(g_1)=p\circ h(g_1)=p\circ f(g_1)=p(g_0)=0.
$$
Since $\varphi$ is a monomorphism, from this we conclude that $p(g_1)=0$; that is, $g_1\in C\times \{0\}=C_0$.
Since $g_1\in G$, it follows that $g_1\in G\cap C_0=G_0$, and so $f_0(g_1)=f(g_1)=g_0$.
This establishes the surjectivity of $f_0$.
$\Box$

\medskip

From \eqref{eq:7d}, we obtain that $h_0=h\restriction_{C_0}:C_0\to C_0$ is a well-defined map. 

\begin{claim}
\label{claim:2}
$h_0:C_0\to C_0$ is a homeomorphism.
\end{claim}

Proof.
Note that $h(G_0)=f(G_0)=f_0(G_0)=G_0$, where the last inclusion follows from Claim \ref{claim:1}.
Since $G$ is dense in $H$ and $C_0=C\times\{0\}$ is open in $H$, it follows that $G_0= G\cap C_0$
is dense in $C_0$. Since $G_0=h(G_0)\subseteq h(C_0)\subseteq C_0$, we conclude that 
$h(C_0)$ is dense in $C_0$.
Therefore, the homomorphism
$h_0=h\restriction_{C_0}:C_0\to C_0$ has dense image.

Without loss of generality, 
we can identify the topological group $C_0$ with $\R^m$
and $h_0$ with a continuous group homomorphism from $\R^m$ to itself. By Fact \ref{fact:5.1},
$h_0$ is a linear map from $\R^m$ to itself.
This means that the image $h_0(\R^m)$ is a closed linear subspace of $\R^m$. Since this image is also dense in $\R^m$ by the previous paragraph, we conclude that $h_0$ is surjective.
By Fact \ref{fact:linear:surjections},
$h_0$ has an inverse which is also a linear map;
that is, $h_0$ is an automorphism of the vector space $\R^m$.
In particular, $h_0$ is a homeomorphism of $C_0$ onto itself.
$\Box$

\medskip
Let $d\in D$ be arbitrary. Since $h$ is a homomorphism
mapping $C_0=C\times \{0\}$ homeomorphically onto itself by Claim \ref{claim:2},
$h$ also maps 
$C\times \{d\}=(0,d)+C\times \{0\}$  homeomorphically onto $C\times\{\varphi(d)\}$. 

Since $\varphi:D\to D$ is a bijection and the family $\{C\times\{d\}:d\in D\}$ forms a clopen partition of $H$, we conclude that $h$ is a homeomorphism of $H$ onto itself. In particular, the restriction $h\restriction_G=f$ of $h$ to $G$ is a homeomorphism of $G$ onto itself. 

We has proved that every continuous automorphism $f$ of $G$ is a homeomorphism. Therefore, $G$ is $g$-reversible
by Definition 
\ref{def:reversible}.
$\Box$

\medskip
A stronger result holds for {\em closed\/} subgroups of $\R^n$.
(Recall that reversibility is stronger than $g$-reversibility by Proposition \ref{reversible_g_reversible}.)

\begin{proposition}
Each closed subgroup of  $\mathbb R^n$ is reversible.
 \end{proposition} 

Proof.
It follows from Fact \ref{fact:dikranjan} that a closed subgroup of $\R^n$ is either locally Euclidean or discrete, so it is reversible 
by Fact~\ref{reversible:fact}(i),(iii).
$\Box$

\medskip

The word ``closed'' cannot be omitted from 
this proposition,
as $\Q$ is a (dense) subgroup of $\R$ which is not reversible; see 
Fact \ref{reversible:fact}(iv).

\section{Additional examples of non-$g$-reversible groups}
\label{sec:6}

First examples of non-$g$-reversible groups were given in Section~\ref{sec:4}. In this section we give two more examples. In Example ~\ref{separable_metrizable_countable_non-g-reversible} 
we exhibit a countable dense subgroup of the Hilbert space which is not $g$-reversible. 
This example 
is attractive due to its simplicity and an explicit description of the subgroup itself. Another example, whose existence is proved in Theorem~\ref{precompact:countable:non-sm}, adds precompactness to the list of properties of the non-$g$-reversible group
but at the expense of both its simplicity and concreteness.

\begin{example}\label{separable_metrizable_countable_non-g-reversible}
{\em The Hilbert space $l^2$ has a countable dense subgroup $G$ which is not $g$-reversible.\/} 
Indeed,
we claim that
$$
G = \{\overline{x} = (x_i) \in l_2: x_i \in \mathbb Q \mbox{ and } x_i = 0 \mbox{ for all but finitely many } i\}
$$
is such a subgroup. Clearly, 
$G$ is  dense in $l_2$  and countable.
Define a continuous
automorphism
$h : G \to G$ of $G$ 
by $h( \overline{x}) = (x_1, \frac{1}{2} x_2, \frac{1}{3} x_3, \dots)$
for $\overline{x} \in G$. For each $k\geq 1$, let $\overline{x}^k = (x_i)$, where $x_k = 1$ and $x_i = 0$ if $i \ne k$.
Let us note that the sequence $\{h(\overline{x}^k)\}$ converges to the point $\overline{0} = (0,0, \dots) \in G$, while the set $\{\overline{x}^1, \overline{x}^2, \dots\}$ is discrete in $G$. This implies that $h^{-1}$ is not continuous at the point $\overline{0}$. 
Therefore, $G$ is not $g$-reversible by Remark \ref{inverse:map:remark}.
\end{example}

Now we describe a general scheme of constructing non-$g$-reversible precompact groups from non-reversible topological spaces.

Recall (\cite[Exercise 7.1.f]{AT}) that for each Tychonoff space $X$ there exists a free abelian precompact topological group $AP(X)$. The description of this topological group can be found in \cite[Section~1]{Sh2}.

\begin{theorem}
\label{free:precompact}
Let $X$ be a space and $AP(X)$ be the free abelian precompact group of $X$. If $X$ is not reversible, then $AP(X)$ is not $g$-reversible.
\end{theorem}
Proof.
Since $X$ is not reversible, there exists a continuous bijection 
$f:X\to X$ whose inverse $f^{-1}$ is discontinuous.
Let $h:AP(X)\to AP(X)$ be a continuous homomorphism extending $f$.
Then $h$ is an automorphism of $AP(X)$.
However, its inverse $h^{-1}$ is discontinuous, as its restriction
$h^{-1}\restriction_{X}$ to $X$ coincides with the discontinuous map
$f^{-1}$.
By Remark \ref{inverse:map:remark}, $AP(X)$ is not $g$-reversible.
$\Box$

\medskip
For a topological space $X$, we use $w(X)$ and $nw(X)$ to denote the weight and network weight of $X$, respectively. 

Let us recall necessary notions from \cite{Sh1}.
Let $X$ be a set. For a fixed positive integer $n$,
an {\em $n$-ary operation on $X$\/} is a map $\psi$ from a subset of $X^n$
to $X$. When $X$ is a topological space and $\psi$ is continuous with respect to the product topology on $X^n$, then the operation $\psi$ is said to be {\em continuous\/}. 
We shall need the following result from \cite{Sh1}:

\begin{fact}
\label{fact:Shakhmatov}
Assume that $\kappa$ is an infinite cardinal, $\mathscr{T}$ is a topology on a set $X$ and
$\Psi$ is a family of continuous operations on 
$(X,\mathscr{T})$ 
such that $|\Psi|\le\kappa$.
Assume also that $nw(X,\mathscr{T})\le\kappa$,
$\mathscr{E}\subseteq \mathscr{T}$ and $|\mathscr{E}|\le\kappa$.
Then there exists a topology $\mathscr{T}'$ on $X$ such that
$\mathscr{E}\subseteq \mathscr{T}'\subseteq \mathscr{T}$, 
$w(X,\mathscr{T}')\le\kappa$ and all operations from $\Psi$
are continuous on $(X,\mathscr{T}')$.
\end{fact}

With the help of this fact, we can prove the following result.
\begin{theorem}
\label{uplotnenie}
Let $(G,\mathscr{T})$ be a non-$g$-reversible (Hausdorff) group.
Then there exists a (Hausdorff) group topology $\mathscr{T}'$ on $G$ such that $w(G,\mathscr{T}')\le nw(G,\mathscr{T})$
and $(G,\mathscr{T}')$ is not $g$-reversible.   
\end{theorem}

Proof.
Let $\kappa=mw(G,\mathscr{T})$.
Since $(G,\mathscr{T})$ is not $g$-reversible, we can fix a continuous automorphism $h$ of $(G,\mathscr{T})$ which is not open. Therefore, there exists $U\in\mathscr{T}$ such that $h(U)\not\in\mathscr{T}$.

Note that we can view $h$ as a unary ($1$-ary) continuous operation on $(G,\mathscr{T})$. Since $(G,\mathscr{T})$ is a topological group, 
the map $\psi: (G,\mathscr{T})^2\to (G,\mathscr{T})$ defined by 
$\psi(x,y)=xy^{-1}$ for $x,y\in G$, is continuous. This means that 
$\psi$ is a continuous binary ($2$-ary) operation on $(G,\mathscr{T})$. Therefore, $\Psi=\{h,\psi\}$ is a family of continuous operations on $(G,\mathscr{T})$ satisfying 
$|\Psi|=2\le\omega\le\kappa$.

Since 
$\mathscr{T}$ is a Hausdorff topology on $G$, there exists a Hausdorff topology $\mathscr{T}^*$ on $G$ such that $\mathscr{T}^*\subseteq \mathscr{T}$ and $w(G,\mathscr{T}^*)\le nw(G,\mathscr{T})=\kappa$ (cf. \cite[Lemma 3.1.18]{E}). 
 Let $\mathscr{B}$ be a base of the topology $\mathscr{T}^*$ such that $|\mathscr{B}|\le\kappa$. Finally, define $\mathscr{E}=\mathscr{B}\cup\{U\}$.
Since $\mathscr{B}\subseteq \mathscr{T}^*\subseteq \mathscr{T}$ and $U\in\mathscr{T}$, we have $\mathscr{E}\subseteq \mathscr{T}$.
Clearly, $|\mathscr{E}|\le\kappa$, as $\kappa$ is an infinite cardinal.

Applying Fact \ref{fact:Shakhmatov} with $G=X$, we can find a topology $\mathscr{T}'$ on $G$ such that 
$\mathscr{E}\subseteq \mathscr{T}'\subseteq \mathscr{T}$,
$w(G,\mathscr{T}')\le\kappa=nw(G,\mathscr{T})$ and all operations from $\Psi$
are continuous on $(G,\mathscr{T}')$. 
Since $\psi\in \Psi$, this means that $(G,\mathscr{T}')$ is a topological group. Since $h\in\Psi$, this means that $h$ is a continuous automorphism of $(G,\mathscr{T}')$.
Since $\mathscr{B}\subseteq \mathscr{E}\subseteq \mathscr{T}'$
and $\mathscr{B}$ is a base of $\mathscr{T}^*$, it follows that 
$\mathscr{T}^*\subseteq \mathscr{T}'$. Since $\mathscr{T}^*$ is 
Hausdorff, so is $\mathscr{T}'$.

Finally, since 
$U\in\mathscr{E}\subseteq \mathscr{T}'$, the set $U$ is $\mathscr{T}'$-open.
On the other hand,
$h(U)\not\in\mathscr{T}$ and $\mathscr{T}'\subseteq \mathscr{T}$ implies $h(U)\not\in \mathscr{T'}$; that is,
$h(U)$ is not $\mathscr{T}'$-open. This means 
that the automorphism
$h:(G,\mathscr{T}')\to (G,\mathscr{T}')$ is not open.

We have found a continuous automorphism $h$ of $(G,\mathscr{T}')$ which is not open, so $(G,\mathscr{T}')$ is not $g$-reversible by Definition \ref{def:g-reversible}.
$\Box$.

\begin{theorem}
\label{precompact:countable:non-sm}
A countable metric precompact abelian group need not be $g$-reversible.
\end{theorem}
Proof.
The space $\Q$ of rational numbers is not reversible
\cite[Example 3]{RW}.
By Theorem \ref{free:precompact}, the free abelian precompact group $G=AP(\Q)$ over $\mathbb Q$ is not $g$-reversible.
Note that $G$ is countable, precompact and abelian.
Since $G$ is countable, it has a countable network, and so $nw(G)\le\aleph_0$.

Let $\mathscr{T}$ be the (Hausdorff) topology of $G$.
Applying Theorem \ref{uplotnenie}, we can find a (Hausdorff) group topology $\mathscr{T'}$ on $G$ such that $w(G,\mathscr{T}')\le nw(G,\mathscr{T})\le\aleph_0$ and $(G,\mathscr{T}')$ is not $g$-reversible. Since $w(G,\mathscr{T}')\le \aleph_0$, we conclude that
$(G,\mathscr{T}')$ is metrizable.
Since $\mathscr{T}'\subseteq \mathscr{T}$ and $\mathscr{T}$ is a precompact group topology on $G$, so is the weaker topology $\mathscr{T}'$.
$\Box$

\medskip
Since countable spaces are $\sigma$-compact,
this theorem shows that local compactness in Theorem 
\ref{omega_narrow_theorem} is essential and cannot be omitted,
even in the presence of metrizability.

\begin{example}
\label{compact:metric:with:dense:non-g-rev}
{\em There exists a compact metric abelian group $K$ having a countable dense subgroup $G$ which is not $g$-reversible\/}.
Indeed, let $K$ be the completion of a group $G$ constructed in Theorem \ref{precompact:countable:non-sm}.
Since $G$ is precompact, $K$ is compact. Note that
$w(K)=w(G)\le\aleph_0$ (cf. \cite[Fact 1.1]{Sh2}).
Thus, $K$ is metrizable.
\end{example}

\begin{corollary}
A (countable) dense subgroup of a (metric) reversible group need not be even $g$-reversible.
\end{corollary}

This corollary shows that Theorem~\ref{theorem_all_subgroups_R_n} 
does not merely follow from the reversibility of $\R^n$
mentioned in Fact \ref{reversible:fact}(iii).

\section{Hereditary $g$-reversibility}
\label{sec:7}

A topological space $X$ is called {\it hereditarily reversible} if 
all 
its subspaces are reversible \cite{RW}.
An analogue of 
this notion
for topological groups is given in the next definition.

\begin{definition}
\label{def:her:g-rev}
We say that a topological group $G$ is  {\em hereditarily $g$-reversible\/} if every subgroup of $G$ is $g$-reversible (in the subspace topology).
\end{definition}

As one may reasonably expect, this notion is much stronger than
mere $g$-reversibility. Indeed, 
Examples~\ref{separable_metrizable_countable_non-g-reversible} and~\ref{compact:metric:with:dense:non-g-rev} can be restated as follows.
\begin{proposition}
\label{non-her-reversible}
\begin{itemize}
\item[(i)] The additive group of the Hilbert space $l_2$ is $g$-reversible but is not hereditarily $g$-reversible.
\item[(ii)] 
A compact metric abelian group 
need not be hereditarily $g$-reversible.
\end{itemize}
\end{proposition}
Proof.
(i)
The Hilbert space $l_2$ is a Polish group, so it is $g$-reversible
by Corollary \ref{Polish_groups}, yet it is not hereditarily reversible by Example~\ref{separable_metrizable_countable_non-g-reversible}
and Definition~\ref{def:her:g-rev}.

(ii)
The topological group $K$ from Example~\ref{compact:metric:with:dense:non-g-rev} is not hereditarily $g$-reversible by Definition~\ref{def:her:g-rev}.
$\Box$

\medskip
Note that the compact group from item (ii) of Proposition~\ref{non-her-reversible} is reversible by  Fact \ref{reversible:fact}(ii), so also $g$-reversible by Proposition \ref{reversible_g_reversible}.

In \cite[Corollary 2.1]{CHH} the authors observed that an infinite
first countable Hausdorff space
is {\em not\/} hereditarily reversible if and only if it contains
a copy of the space $N_{\aleph_0}$ defined in Fact \ref{reversible:fact}(v). This implies that the non-trivial convergent sequence (together with its limit) and
discrete spaces are
the only hereditarily reversible Hausdorff first countable spaces. 
This scarcity of first countable Hausdorff hereditary reversible spaces is in a sharp contrast with 
the following result which is an immediate consequence of
Theorem~\ref{theorem_all_subgroups_R_n} and
Definition \ref{def:her:g-rev}.
\begin{proposition}
(Every subgroup of) $\R^n$ is hereditarily $g$-reversible. $\Box$
\end{proposition}

This proposition suggests the following question:

\begin{question}
Is the topological group $\mathbb{R}^{\aleph_0}$ hereditarily $g$-reversible? 
\end{question}

Note that by
Anderson's theorem (cf. \cite{vM}),
$\mathbb R^{\aleph_0}$ is homeomorphic to the Hilbert space $l^2$, and the latter is not hereditarily $g$-reversible by 
Proposition~\ref{non-her-reversible}(i).

Our next proposition strengthens the conclusion of Theorem \ref{Bohr_topology}.

\begin{proposition}
Every abelian group equipped with its Bohr topology is hereditarily $g$-reversible.
\end{proposition}

Proof. 
Let $G$ be an abelian group with its Bohr topology and let $H$ be a subgroup of $G$. It is
well known
that the subspace topology of $H$ coincides with the Bohr topology of $H$, so the subgroup $H$ of $G$ is $g$-reversible by
Theorem \ref{Bohr_topology}.
Thus, $G$ is $g$-reversible by Definition~\ref{def:her:g-rev}.
$\Box$

\medskip
Since an infinite precompact group cannot be discrete,
we obtain the following
                                                                                                                        \begin{corollary}
                                                                                                                        Every infinite abelian group admits a non-discrete hereditarily $g$-reversible group topology.
                                                                                                                       \end{corollary}
                                                                                                                          
A topological group $G$ is called {\em hereditarily minimal\/} if every subgroup of $G$ is minimal in its subspace topology \cite{XDST}. Combining this definition with Theorem~\ref{minimal:theorem}, we get the following
\begin{proposition}
\label{her:minimal:are:her:g-reversible}
Hereditarily minimal groups are hereditarily $g$-reversible.
\end{proposition}

Prodanov proved that
an infinite compact abelian group is hereditarily minimal if and only if it is topologically isomorphic to the group $\Z_p$ of $p$-adic integers, 
for some prime number $p$ \cite{Prodanov}. In \cite{DS}, 
Dikranjan and Stoyanov characterized all 
hereditarily minimal abelian topological groups (see also \cite[Fact 1.4]{XDST}).
It follows from 
this characterization that the groups of $p$-adic integers $\Z_p$ (for a prime number $p$) are the only infinite locally compact, hereditarily minimal abelian groups \cite[Corollary 1.5]{XDST}.
The authors of \cite{XDST} made a progress in the direction of describing non-abelian hereditarily minimal topological 
groups. 

The discussion above and Proposition \ref{her:minimal:are:her:g-reversible} 
justify the following two problems.

\begin{problem}
Describe all hereditarily $g$-reversible (abelian) groups.
\end{problem}

The following particular versions of this general problem may be more tractable.

\begin{problem}
Describe hereditarily $g$-reversible 
topological groups in the following classes:
\begin{itemize}
\item[(i)] compact (abelian) groups;
\item[(ii)] locally compact (abelian) groups;
\item[(iii)] separable metric (abelian) groups.
\end{itemize}
\end{problem}

\section{$g$-reversibility in closed subgroups}
\label{sec7:b}

In 
this
section, we turn our attention on a question as to which extent $g$-reversibility of a topological group is inherited by its closed subgroups. 

Recall that a subgroup $H$ of a topological group $G$ is called a {\em group retract of $G$\/} provided that there exists a continuous homomorphism $h: G\to H$ (called a {\em group retraction of $G$ into $H$\/}) satisfying $h(x)=x$ for every $x\in H$.
Note that if $h: G\to H$ is a group retraction, then $h$ is a quotient map and 
$H$ is closed in $G$, so $H$ is both a closed subgroup of $G$ and a quotient group of $G$.

The following fundamental result is due to Megrelishvili \cite[Theorem 7.2]{Me2}.

\begin{fact}
\label{Megrelishvili}
Every topological group $H$ is a group retract of some minimal group $G$ such that $w(G)=w(H)$.
\end{fact}

Since minimal groups are $g$-reversible by Theorem~\ref{minimal:theorem}, from this fact we get the following

\begin{corollary}
\label{closed:subgroups}
Every topological group $H$ is a group retract of some $g$-reversible group $G$ such that $w(G)=w(H)$.
In particular, $H$ is both a closed subgroup of $G$ and a quotient group of $G$.
\end{corollary}

\begin{corollary}
\label{cor:7.8:a}
A (countable precompact) closed subgroup of a $g$-reversible (separable metric) group need not be $g$-reversible.
\end{corollary}

Proof. Let $H$ 
be a countable metric precompact non-$g$-reversible abelian group constructed in 
Theorem~\ref{precompact:countable:non-sm}. Applying Corollary~\ref{closed:subgroups} to this $H$, we get a  $g$-reversible group $G$ containing $H$ as a closed subgroup such that $w(G)=w(H)=\aleph_0$. Since $G$ has a countable base, it is separable metric.
$\Box$

\begin{example}
\label{loc:comp:ex:of:closed:subgroup}
{\em There exists a locally compact metric minimal (thus, $g$-reversible by Theorem~\ref{minimal:theorem}) topological group which contains a closed non-$g$-reversible subgroup.\/}
Indeed, let $H$ be a non-$g$-reversible locally compact metric abelian group from Corollary~\ref{corollary:4.4}. It was proved in \cite{Me1} that $H$ is a group retract of a locally compact metric 
minimal (thus, $g$-reversible) group $G$. In particular, $H$ is closed in $G$.
\end{example}

The ambient groups in Fact~\ref{Megrelishvili}, Corollaries~\ref{closed:subgroups}, \ref{cor:7.8:a} and Example~\ref{loc:comp:ex:of:closed:subgroup} are non-abelian.
This leaves open the following pair of questions:

\begin{question}
\label{closed:subgroup:question}
\label{closed:subgroup:question:lc}
\begin{itemize}
\item[(i)] Must a closed subgroup 
of a $g$-reversible 
abelian group $G$
be $g$-reversible? What is the answer if $G$ is assumed to be even reversible? 
\item[(ii)]
Must every closed subgroup of a $g$-reversible locally compact abelian group $G$ be $g$-reversible? What is the answer if $G$ is assumed to be even reversible? 
\end{itemize}
\end{question}

\begin{remark}
{\em A closed subgroup of a 
$g$-reversible
topological group $G$ is 
$g$-reversible
in each of the following cases:
\begin{itemize}
\item[(i)] $G$ is a Polish group;
\item[(ii)] $G$ is a $\sigma$-compact locally compact group.
\end{itemize}
}
Indeed, (i) follows from Corollary \ref{Polish_groups} and the fact that a closed subgroup of a Polish group is Polish.
Item (ii) follows from Theorem \ref{omega_narrow_theorem} and the fact that a closed subgroup of a $\sigma$-compact locally compact group is itself $\sigma$-compact and locally compact.
\end{remark}

\section{$g$-reversibility and topological products}
\label{sec:8}

In this section we consider a question of preservation of $g$-reversibility under taking Tychonoff products. 

We start with the necessary condition.

\begin{theorem}
\label{products:necessary}
Let $G=\prod_{i\in I} G_i$ be the Tychonoff product of a family $\{G_i:i\in I\}$ of topological groups. If $G$ is $g$-reversible, 
then every factor $G_i$ must be $g$-reversible as well.
\end{theorem}

Proof. Fix $i\in I$ and a continuous automorphism $f_i$ of $G_i$.
For $j\in I\setminus\{i\}$, let $f_j$ be the identity automorphism of $G_j$. Then the map $f: G\to G$ defined by $f(x)=(f_i(x_i))$ for 
every $x=(x_i)\in G$, is a continuous automorphism of $G$.
Since $G$ is $g$-reversible, the map $f$ is open. Now one can easily conclude from this and the definition of $f$ that $f_i$ is open as well.
This shows that $G_i$ is $g$-reversible.
$\Box$

\medskip
The following corollary allows one to produce new non-$g$-reversible groups from a given one.

\begin{corollary}\label{non-g-reversible_groups} Let $G_1$ be a non-$g$-reversible topological group and $G_2$ be any topological group. Then the product $G_1 \times G_2$ is not $g$-reversible.
\end{corollary}

\begin{remark}
The converse to Theorem \ref{products:necessary} 
does not hold.
Since both discrete groups and compact groups are reversible
by Fact \ref{reversible:fact}(i),(ii), 
Corollary~\ref{g-nonreversible_product} shows that the topological product of two (locally compact) reversible abelian groups need not be $g$-reversible.
It remains only to recall that  reversible groups are $g$-reversible by Proposition \ref{reversible_g_reversible}.
\end{remark}

\begin{proposition}\label{product_theorem} Let $G_1$ and $G_2$ be $\sigma$-compact locally compact topological groups.
Then the topological product $G = G_1 \times G_2$ is $g$-reversible.
\end{proposition}
Proof. Note that the topological group $G$ is $\sigma$-compact and locally compact, so we can apply Theorem~\ref{omega_narrow_theorem}. $\Box$

\medskip
Corollary~\ref{g-nonreversible_product}  shows that
the assumption of 
$\sigma$-compactness in Proposition~\ref{product_theorem} is essential.

\begin{question} Let $G_1$ and $G_2$ be connected $g$-reversible (or even reversible) topological groups. Is the topological product
$G_1 \times G_2$ $g$-reversible?
\end{question}

Even the quite special case of topological products of a single discrete group and a single compact group appears to be interesting. 

\begin{problem}
Describe the class $\mathbf{K}$ of cardinals $\kappa$ such that
all topological products $D\times K$, where $D$ is a discrete group of cardinality $\kappa$ and $K$ is a compact group, are $g$-reversible. 
\end{problem}

Clearly, $\omega\subseteq \mathbf{K}$, as the topological product of a finite group with a compact group is compact, and so $g$-reversible by 
Remark~\ref{rem:2.10}. Furthermore, $\omega\in\mathbf{K}$, as 
the topological product $G$ of a countable discrete group with a compact group is 
locally compact and $\sigma$-compact; now the $g$-reversibility of $G$ follows from Theorem~ \ref{omega_narrow_theorem}.
Corollary~\ref{g-nonreversible_product} shows that $\mathfrak{c}\not\in \mathbf{K}$. Therefore, $\omega_1\not\in\mathbf{K}$ under the assumption of the Continuum Hypothesis CH.

\begin{question}
Can $\omega_1\not\in\mathbf{K}$ be proved in ZFC?
\end{question}

\begin{question}
Does the equality $\mathbf{K}=\omega\cup\{\omega\}$ hold? In other words, can one find, for every uncountable cardinal $\kappa$, a discrete group $D$ of cardinality $\kappa$ such that its product 
$D\times K$ with  some compact group $K$ is not $g$-reversible?
\end{question}

\begin{question}
Let $\{G_i:i\in I\}$ be an arbitrary family of Polish groups.
Is the product $\prod_{i\in I} G_i$ $g$-reversible?
\end{question}

Since countable products of Polish groups are Polish, Corollary~\ref{Polish_groups} provides a positive answer 
to this question for an at most countable set $I$.

Even the following concrete version of this question remain unclear.
\begin{question}
Are topological groups $\Z^\kappa$ and $\R^\kappa$ $g$-reversible for an arbitrary (uncountable) cardinal $\kappa$?
\end{question}

\begin{remark}
Let $\kappa$ be a cardinal 
less than
the first measurable cardinal. Then every homomorphism from $\Z^\kappa$ to $\Z$ 
is continuous; see discussion in~\cite[Sec. 2]{Nunke}.
It easily follows from this that every automorphism of $\Z^\kappa$ 
is continuous, and so $\Z^\kappa$ is $g$-reversible.
Therefore, {\em if there are no measurable cardinals, then 
$\Z^\kappa$ is $g$-reversible for every cardinal $\kappa$.\/}
\end{remark}

\begin{question}
For every cardinal $\kappa\ge 1$, is there a topological (abelian) group $G$, depending on $\kappa$, such that $G^\lambda$ is $g$-reversible for every cardinal $\lambda<\kappa$, yet $G^{\kappa}$ is not $g$-reversible?  
\end{question}

\section{Comparison of reversible spaces with $g$-reversible groups}
\label{sec:9}

In this short section we give four examples contrasting reversibility with $g$-reversibility in particular topological groups. 

\begin{example}
\label{ex:9.1}
Recall that no infinite-dimensional normed space is reversible, see \cite[Theorem 2]{RW}. In particular, the Hilbert space $l_2$ is not reversible. On the other hand, $l_2$ is a Polish group, so it is $g$-reversible by Corollary \ref{Polish_groups}.
\end{example}

\begin{example}
\label{example:9.2}
  {\em Let $G$ be a subgroup of $\R^n$. Then the topological group $\Q \times G$ is $g$-reversible but is not reversible as a topological space.\/}
Indeed,
the topological group $\Q \times G$ is topologically isomorphic to a  subgroup of $\R^{n+1}$, so it is  $g$-reversible by Theorem~\ref{theorem_all_subgroups_R_n}. By Example~\ref{the group of rational numbers}, the topological space $\Q$ is non-reversible.
Hence, the Cartesian product $\Q \times G$ is also non-reversible
by Fact~\ref{reversible:fact}(viii).
\end{example} 

We now compare two notions in products $\Z^{\kappa} \times \{0,1\}^{\lambda} \times \R^\mu$, for cardinals $\kappa,\lambda,\mu$.

\begin{example}
Assume that $\max\{\kappa,\lambda,\mu\}\le\aleph_0$.
Then the topological group $G=\Z^{\kappa} \times \{0,1\}^{\lambda} \times \R^\mu$
is Polish, so it is
$g$-reversible by Corollary~ \ref{Polish_groups}.
Next, we discuss the reversibility of $G$. 

If all three cardinals $\kappa,\lambda,\mu$ are finite, then $G$ is locally Euclidean, so reversible by Fact~\ref{reversible:fact}(iii).

If $\kappa = \aleph_0$,  then 
$G$ contains the factor $\mathbb Z^{\aleph_0}$ which is homeomorphic to the space $\mathbb P$ of irrational numbers
which is not reversible by \cite[Example 4]{RW}, so $G$ is not reversible as a topological space
by Fact~\ref{reversible:fact}(viii). 

If $\mu = \aleph_0$,  then 
$G$
contains the factor $\mathbb R^{\aleph_0}$ which is homeomorphic to the space $l^2$ by 
Anderson's theorem (cf. \cite{vM}). Since $l_2$ is not reversible by Example~\ref{ex:9.1},
$G$ is not reversible 
by Fact~\ref{reversible:fact}(viii). 

If $\lambda = \aleph_0$ and $\kappa \geq 1$, then 
$G$
contains the factor $\mathbb Z \times \mathbb C$, where $\mathbb C$ is the Cantor set, which is not reversible, so
$G$
is not reversible
by Fact~\ref{reversible:fact}(viii). 

If $\lambda = \aleph_0$, $\kappa =0$ and $\mu = 0$, then $G$ is compact and so is reversible by  Fact~\ref{reversible:fact}(ii).

When
$\lambda = \aleph_0$, $\kappa =0$ and $\mu > 0$, the reversibility of $G$ is unknown.
\end{example}

\begin{example} 
Let $n,m$ be non-negative integers and $\lambda$ be a cardinal.
Then 
the topological group $G=\Z^m \times \{0,1\}^{\lambda} \times \R^n$
is $\sigma$-compact and locally compact, so it is $g$-reversible
by Theorem~\ref{omega_narrow_theorem}.

On the other hand, if $m\ge 1$ and $\lambda\ge\aleph_0$, 
then $G$ is not reversible.
Indeed, $G$ contains as a factor the product $\Z\times\{0,1\}^{\aleph_0}$, and since
$\{0,1\}^{\aleph_0}$ is homeomorphic to the Cantor set $\mathbb C$, $G$ contains as a factor the topological product $\Z \times \mathbb C$. Since the latter is not reversible, 
$G$ is not 
reversible as a topological space by Fact~\ref{reversible:fact}(viii). 
\end{example} 

\section{Concluding remarks and additional open problems}
\label{sec:10}

Recall that an abelian group $G$ is called {\em rigid\/} if $\id_G$ and $-\id_G$ are the only automorphisms of $G$. 

The following remark is due to Dekui Peng.

\begin{remark}
\label{Peng:remark}
(i) {\em All 
group topologies on a  rigid abelian group $G$ are $g$-reversible\/}. 
Indeed,
both $\id_G$ and $-\id_G$ are homeomorphisms with respect to any group topology on $G$. 

(ii) The group $\Z$ of integer numbers is a well-known example of a rigid abelian group. Therefore, $\Z$ is an example of an abelian group 
on which all group topologies are $g$-reversible. 
\end{remark}

\begin{remark} 
\label{Q:remarks}
Note that the topological group $G$ from  Example~\ref{separable_metrizable_countable_non-g-reversible} is homeomorphic to the space $\mathbb Q$ of rational numbers. Together with Example~\ref{the group of rational numbers}  this shows that {\em the space $\mathbb{Q}$ of rational numbers admits two group structures such that one topological group is $g$-reversible and the other one is not\/}.
By Example~\ref{example:9.2}
and Corollary~\ref{non-g-reversible_groups},  the topological products $\mathbb Q \times \mathbb R^n, n \geq 1,$ possess the same property.
\end{remark}

\begin{question}What topological spaces admit two group structures such that one topological group is $g$-reversible and the other one is not?
\end{question}

\begin{problem}
\begin{itemize}
\item[(i)]
Describe $g$-reversible pseudocompact (abelian) groups.
\item[(ii)] Describe $g$-reversible countably compact (abelian) groups.
\end{itemize}
\end{problem}

A concrete version of this problem asks whether the converse of Theorem~\ref{minimal:theorem} holds for pseudocompact 
or countably compact
(abelian) groups. 

\begin{question}
Must every pseudocompact (countably compact) $g$-reversible group $G$ be minimal? What is the answer to  this question if one additionally assumes that $G$ is abelian?
\end{question}

\begin{remark}
It seems interesting to consider analogues of reversibility
in other categories, for example, topological semigroups, topological rings, topological vector spaces etc. The introduction of these new notions naturally follows the similar path from 
Definition~\ref{def:reversible} to Definition~\ref{def:g-reversible}
which we executed in this paper. One may expect the results for other categories would not only be similar to those for the category of topological groups but also reflect some particularities of these categories.
\end{remark}

\medskip

\noindent
{\bf Acknowledgements:\/}
We are grateful to Dekui Peng
for his Remark~\ref{Peng:remark}.
We are obliged to Takamitsu Yamauchi for providing a stimulating reference, 
as well as to Menahem Shlossberg for pointing out reference \cite{XDST}.
At last but not least,  
the first listed author would like to thank 
Alexander Karassev
for his inspiration of Example~\ref{separable_metrizable_countable_non-g-reversible}.

\vskip1cm

\noindent(V.A. Chatyrko)\\
Department of Mathematics, Linkoping University, 581 83 Linkoping, Sweden.\\
vitalij.tjatyrko@liu.se

\vskip0.3cm
\noindent(D.B. Shakhmatov) \\
Graduate School of Science and Engineering, 
Division of Science,
Ehime University\\
Bunkyo-cho 2-5,
790-8577 Matsuyama,
Japan\\
dmitri.shakhmatov@ehime-u.ac.jp

\end{document}